\documentclass[a4paper,10pt]{article}
\usepackage{amsmath,amsthm,xypic,
amssymb,amsfonts,
pdfsync,stmaryrd,mathrsfs,yfonts}
\usepackage[capitalise]{cleveref}
\newtheorem{thm}{Theorem}[section]

\newtheorem{prop}[thm]{Proposition}

\newtheorem{lem}[thm]{Lemma}
\theoremstyle{definition}

\newtheorem{rem}[thm]{Remark}

\theoremstyle{plain}

\newcommand{\gal}{\mathop{\mathrm{Gal}}}
\newcommand{\co}{\mathop{\mathrm{cor}}}

\renewcommand{\phi}{\varphi}

\newcommand{\symd}{\mathop{\mathrm{Symd}}}
\newcommand{\sym}{\mathop{\mathrm{Sym}}}

\newcommand{\Hom}{\mathop{\mathrm{Hom}}}

\newcommand{\car}{\mathop{\mathrm{char}}}
\newcommand{\id}{\mathop{\mathrm{id}}}

\newcommand{\an}{\mathop{\mathrm{an}}}

\author{A.-H. Nokhodkar}
\title
{Quadratic descent of hermitian forms}

\begin{document}
\date{}
\maketitle
\begin{abstract}
\noindent
  Quadratic descent of hermitian and skew hermitian forms over division algebras with involution of the first kind in arbitrary characteristic is investigated and a criterion, in terms of systems of quadratic forms, is obtained.
A refined result is also obtained for hermitian (resp. skew hermitian) forms over a quaternion algebra with symplectic (resp. orthogonal) involution.
\\
\\
\noindent
\emph{Mathematics Subject Classification:} 11E39, 11E04. \\
\emph{Keywords:}  Hermitian form, quadratic descent, system of quadratic forms, quaternion algebra.\\
\end{abstract}
\section{Introduction}
The classical Springer theorem asserts that if $K/F$ is a field extension of odd
degree, then the canonical homomorphism of Witt rings $r_{K/F}:W(F)\rightarrow W(K)$ is injective.
A theorem proved by Rosenberg and Ware \cite{rosen} states that if $K/F$ is a Galois extension of odd
degree in characteristic not two, then $r_{K/F}:W(F)\rightarrow W(K)^{\gal(K/F)}$
is an isomorphism.
The image of $r_{K/F}$ was studied further in \cite{rost} for extensions of odd degree in characteristic not two and a descent property of Pfister forms was obtained.
Also, for purely inseparable extensions, some descent properties of the associated Witt
group of hermitian forms of a central simple algebra with involution over $K$ were studied in \cite{bayer}.

In this work, we investigate quadratic descent of hermitian and skew hermitian forms over division algebras with involution of the first kind in arbitrary characteristic.
Our study is based on some constructions of \cite{meh}.
For $\lambda=\pm1$, let $h$ be a $\lambda$-hermitian form over a division algebra with involution of the first kind $(D,\theta)$.
For every basis $\mathcal{B}$ of the set of $\lambda$-symmetrized elements of $(D,\theta)$, a system of quadratic forms $q_{h,\mathcal{B}}$ was associated to $h$ in \cite{meh}.
This system is a generalization of the so-called Jacobson's trace form defined in \cite{jacob} (see also \cite{sah}).
It was shown that, except for certain special cases, the system $q_{h,\mathcal{B}}$ determines the isotropy and metabolicity of $h$, as well as the isometry class of $h$ (see \cite[\S 5]{meh}).
In view of these characterizing properties, a natural problem is to find a criterion for a hermitian form $h$, in terms of $q_{h,\mathcal{B}}$, to have a quadratic descent.
Our main result is Theorem \ref{main} which states that for a separable quadratic extension $K/F$, a division $F$-algebra with involution of the first kind $(D,\theta)$ and an even $\lambda$-hermitian space $(V,h)$ over $(D,\theta)_K$ the following statements are equivalent (except for the case where $D=F$ and $\lambda=-1$):
\begin{itemize}
  \item [$(i)$] $h$ has a descent to $(D,\theta)$.
  \item [$(ii)$] There exists a basis $\mathcal{B}$ of $\symd_\lambda(D,\theta)$ for which $q_{h,\mathcal{B}}$ has a descent to $F$.
  \end{itemize}

In \S5 we study in more detail quadratic descent of $\lambda$-hermitian forms over quaternion division algebras with involution of the first kind whose corresponding systems of quadratic forms reduce to quadratic forms.
For a field extension $K/F$ and a quaternion division algebra with involution of the first kind $(Q,\sigma)$ over $K$,
we will find in Theorem \ref{main2} a criterion for the aforementioned hermitian forms over $(Q,\sigma)$ to have descents to $F$.
This criterion is stronger than Theorem \ref{main}, in the sense that the underlying algebra with involution $(Q,\sigma)$ is now defined over $K$, while  the pair $(D,\theta)_K$ in Theorem \ref{main} was extended from $(D,\theta)$ over $F$.

\section{Preliminaries}
Let $A$ be a central simple algebra over a field $F$.
An {\it involution} on $A$ is an antiautomorphism of $A$ of order two.
An involution $\sigma$ on $A$ is said to be of {\it the first kind} if it restricts to the identity on $F$.
Otherwise, it is said to be of the {\it second kind}.
For a central simple algebra with involution $(A,\sigma)$ and $\lambda=\pm1$ we use the notation
\begin{align*}
{\sym}_\lambda(A,\sigma)&=\{x\in A\mid \sigma(x)=\lambda x\},\\
{\symd}_\lambda(A,\sigma)&=\{x+\lambda\sigma(x)\mid x\in A\}.
\end{align*}
Note that ${\sym}_\lambda(A,\sigma)\subseteq {\symd}_\lambda(A,\sigma)$ with equality if $\car F\neq 2$.

Let $K/F$ be a field extension and let $(A,\sigma)$ be a central simple algebra with involution over $K$.
We say that $(A,\sigma)$ has a {\it descent} to $F$ if there exists a central simple algebra with involution $(A',\sigma')$ over $F$ such that $(A,\sigma)\simeq(A',\sigma')_K$, where $(A',\sigma')_K:=(A'_K,\sigma'_K)$ is the $K$-algebra with involution obtained from $(A',\sigma')$ by extending scalars to $K$, i.e., $A'_K=A'\otimes_FK$ and $\sigma'_K\simeq \sigma'\otimes\id$.

Let $(A,\sigma)$ be a central simple algebra with involution of the first kind over a field $F$.
According to \cite[(2.1)]{knus} for every splitting field $L$ of $A$, the involution $\sigma_L$ becomes adjoint to a symmetric or antisymmetric bilinear form $\mathfrak{b}$ over $L$.
We say that $\sigma$ is {\it symplectic} if $\mathfrak{b}$ is alternating and {\it orthogonal} otherwise.
Also, the involution $\sigma$ is said to be {\it of type} $1$ (resp. {\it of type} $-1$) if $\mathfrak{b}$ is symmetric (resp. antisymmetric).
Note that  an involution of the first kind over a field of characteristic different from $2$ is of type $1$ if and only if it is orthogonal.

Let $(D,\theta)$ be a finite dimensional division algebra with involution over a field $F$ and let $\lambda=\pm1$.
Let $V$ be a finite dimensional right vector space over $D$.
A {\it $\lambda$-hermitian form} on $V$ is a bi-additive map $h:V\times V\rightarrow D$ satisfying $h(u\alpha,v\beta)=\theta(\alpha) h(u,v)\beta$ and $h(v,u)=\lambda\theta(h(u,v))$ for all $u,v\in V$ and $\alpha,\beta\in D$.
The pair $(V,h)$ is called a {\it $\lambda$-hermitian space} over $(D,\theta)$.
If $\lambda=1$ (resp. $\lambda=-1$), $h$ is also called a {\it hermitian} (resp. {\it skew hermitian}) {\it form}.
For a subspace $W$ of $V$, the {\it orthogonal complement} of $W$ (with respect to $h$) is defined as
\[W^{\perp_h}=\{v\in V\mid h(v,w)=0 \ {\rm for\ all}\ w\in W \}.\]
The form $h$ is called {\it regular} if $V^{\perp_h}=\{0\}$.
A $\lambda$-hermitian form $h$ over $(D,\theta)$ is called {\it even} if $h(v,v)\in\symd_\lambda(D,\theta)$ for all $v\in V$.

Let $(D,\theta)$ be a division algebra with involution over $F$ and let $(V,h)$ be a $\lambda$-hermitian space over $(D,\theta)$.
If $K/F$ is a finite extension for which $D_K$ is a division algebra then there exists a $\lambda$-hermitian space $(V_K,h_K)$ over $(D,\theta)_K$, where $V_K=V\otimes_FK$ and $h_K:V_K\times V_K\rightarrow D_K$ is induced
by $h_K(u\otimes\alpha,v\otimes\beta)=h(u,v)\otimes\alpha\beta$ for $u,v\in V$ and $\alpha,\beta\in K$.

Let $K/F$ be a field extension and let $(D,\theta)$ be a division algebra with involution of the first kind over $K$.
For $\lambda=\pm1$, let $(V,h)$ be a $\lambda$-hermitian space over $(D,\theta)$.
We say that $h$ {\it has a descent to} $F$ if there exist a division algebra with involution $(D',\theta')$ over $F$ and a $\lambda$-hermitian space $(V',h')$ over $(D',\theta')$ such that $(D,\theta)\simeq(D',\theta')_K$ and $(V,h)\simeq (V'_K,h'_K)$.
Further, if $(D'',\theta'')$ is a division algebra with involution over $F$ satisfying $(D,\theta)\simeq(D'',\theta'')_K$, we say that $h$ {\it has a descent to} $(D'',\theta'')$ if there exists a $\lambda$-hermitian space $(V'',h'')$ over $(D'',\theta'')$ such that $(V,h)\simeq (V''_K,h''_K)$.

\section{Quadratic maps}
Let $V$ be a finite dimensional vector space over a field $F$ of arbitrary characteristic.
A {\it quadratic form} on $V$ is a map $q:V\rightarrow F$ for which (i) $q(\alpha v)=\alpha^2q(v)$ for all $\alpha\in F$ and $v\in V$;
(ii) the map $\mathfrak{b}_q:V\times V\rightarrow F$ given by $\mathfrak{b}_q(u,v)=q(u+v)-q(u)-q(v)$ is a bilinear form.
The pair $(V,q)$ is called a {\it quadratic space} over $F$.
The {\it orthogonal complement} of a subspace $W$ of $V$ is defined as
\[W^\perp=\{v\in V\mid \mathfrak{b}_q(v,w)=0\ {\rm for\ all}\ w\in W\}.\]
We say that $q$ is {\it regular} if $V^\perp=\{0\}$.
The form $q$ is called {\it isotropic} if there exists a nonzero  vector $v\in V$ such that $q(v)=0$ and {\it anisotropic} otherwise.
Let $K/F$ be a field extension and let $(V,q)$ be a quadratic space over $F$.
We say that $(V,q)$ (or the quadratic form $q$  itself) has a {\it descent} to $F$ if there exists a quadratic space $(V',\rho)$ over $F$ such that $(V,q)\simeq(V'_K,\rho_K)$, where $\rho_K:V'_K\rightarrow K$ is the quadratic form satisfying $\rho_K(v\otimes\alpha)=\alpha^2\rho(v)$ for all $\alpha\in K$ and $v\in V'$.

An {\it $m$-fold system of quadratic forms} on $V$ is an $m$-tuple $q=(q_1,\cdots,q_m)$, where every $q_i:V\rightarrow F$ is a quadratic form.
The system $q$ may be identified with a quadratic map $q:V\rightarrow F^m$ (see \cite[p. 132]{pfister}).
This quadratic map induces a polar map $\mathfrak{b}_q:V\times V\rightarrow F^m$ defined by
\[\mathfrak{b}_q(u,v)=q(u+v)-q(u)-q(v).\]
The system $q$ is called {\it regular} if there is no nonzero vector $v\in V$ such that $\mathfrak{b}_q(u,v)=0\in F^m$ for all $u\in V$.
We say that $q$ is {\it totally regular} if every $q_i$ is regular.
Two systems of quadratic forms $q:V\rightarrow F^m$ and $q':V'\rightarrow F^m$ are said to be {\it equivalent} if there exists a linear isomorphism $f:V\rightarrow V'$ such that $q'(f(v))=q(v)$ for every $v\in V$.
In this case, we write $(V,q)\simeq (V',q')$ or simply $q\simeq q'$.

Let $K/F$ be a field extension.
Let $V$ be a finite dimensional vector space over $K$ and let $q:V\rightarrow K^m$ be a quadratic map.
We say that $q$ has a {\it descent} to $F^m$ if there exist a vector space $V'$ over $F$ and a quadratic map $\rho:V'\rightarrow F^m$ such that $q\simeq\rho_K$,
where $\rho_K:V'_K\rightarrow K^m$ is the quadratic map induced by $\rho_K(v\otimes\alpha)=\rho(v)\alpha$ for $v\in V'$ and $\alpha\in K$.
Note that if $q=(q_1,\cdots,q_m)$ then $q$ has a descent to $F^m$ if and only if every $q_i$ has a descent to $F$.

Let $K/F$ be a finite field extension and let $s:K\rightarrow F$ be an $F$-linear functional.
If $(V,\mathfrak{b})$ is a bilinear space over $K$ the {\it transfer} $s_*(\mathfrak{b})$ of $\mathfrak{b}$ is a symmetric bilinear form on $V$, as a vector space over $F$, defined by
$s_*(\mathfrak{b})(u,v)=s(\mathfrak{b}(u,v))$ for $u,v\in V$.
Also,  if $q$ is a quadratic form on $V$, the {\it transfer} $s_*(q)$ of $q$  is a quadratic $F$-form on $V$ whose polar form is $s_*(\mathfrak{b_q})$  and satisfies
$(s_*(q))(v)=s(q(v))$ for $v\in V$.

\begin{lem}\label{lem}
Let $K/F$ be a separable quadratic extension and let $s:K\rightarrow F$ be a nonzero $F$-linear functional with $s(1)=0$.
A regular quadratic form $q$ over $K$ has a descent to $F$ if and only if $s_*(q)$ is hyperbolic.
\end{lem}

\begin{proof}
Observe first that if $\car F=2$ then $q$ is even dimensional.
Hence, the `only if' implication follows from \cite[(34.4)]{elman} if $\car F\neq2$ and the exactness of the sequence in \cite[(34.9)]{elman} at $I_q(K)$ if $\car F=2$.
To prove the converse, write $q\simeq q_{\an}\perp k\mathbb{H}$, where $q_{\an}$ is anisotropic, $\mathbb{H}$ is the hyperbolic plane and $k$ is a nonnegative integer.
By \cite[(20.1)]{elman}, we have $s_*(q)\simeq s_*(q_{\an})\perp s_*(k\mathbb{H})$.
Also, $s_*(k\mathbb{H})$ is hyperbolic by \cite[(20.5)]{elman}.
Hence, $s_*(q_{\an})$ is also hyperbolic.
The proofs of \cite[(34.4)]{elman} (if $\car F\neq2$) and \cite[(34.9)]{elman} (if $\car F=2$) show that $q_{\an}$ has a descent to $F$.
Hence, $q$ has a descent to $F$, proving the converse implication.
\end{proof}

\begin{prop}\label{desq}
Let $K/F$ be a separable quadratic field extension and let $q\simeq \rho\perp \phi$ be an isometry of totally regular quadratic maps over $K^m$.
If $q$ and $\rho$ have descents to $F^m$, then so does $\phi$.
\end{prop}

\begin{proof}
Write $q=(q_1,\cdots,q_m)$, $\rho=(\rho_1,\cdots,\rho_m)$ and $\phi=(\phi_1,\cdots,\phi_m)$ for some quadratic forms $q_i$, $\rho_i$ and $\phi_i$ over $K$, $i=1,\cdots,m$.
Then  all $q_i$ and $\rho_i$ have descents to $F$.
We want to show that for $i=1,\cdots,m$, $\phi_i$ has a descent to $F$.
Fix an index $i$ with $1\leqslant i\leqslant m$.

Let $s:K\rightarrow F$ be a nonzero $F$-linear functional with $s(1)=0$.
Since $q_i\simeq \rho_i\perp \phi_i$, by \cite[(20.1)]{elman} we have $s_*(q_i)\simeq s_*(\rho_i)\perp s_*(\phi_i)$.
Also, $s_*(q_i)$ and $s_*(\rho_i)$ are hyperbolic by \cite[(34.9)]{elman}, which implies that $s_*(\phi_i)$ is also hyperbolic.
The result now follows from
Lemma \ref{lem}.
\end{proof}

\section{Quadratic descent of hermitian forms}
We recall here some constructions from \cite{meh}.

Let $(D,\theta)$ be a finite dimensional division algebra with involution of the first kind over a field $F$.
Let $(V,h)$ be a $\lambda$-hermitian space over $(D,\theta)$, where $\lambda=\pm1$.
Consider a basis $\mathcal{B}=\{u_1,\cdots,u_n\}$ of $\sym_{\lambda}(D,\theta)$ over $F$ and let $\{\pi_1,\cdots,\pi_n\}$ be its dual basis of $\Hom(\sym_\lambda(D,\theta),F)$.
For $i=1,\cdots,n$, define a map
\[q_{h,\mathcal{B}}^{u_i}:V\rightarrow F\]
via
\[q_{h,\mathcal{B}}^{u_i}(v)=\pi_i(h(v,v)).\]
Let $q_{h,\mathcal{B}}=(q_{h,\mathcal{B}}^{u_1},\cdots,q_{h,\mathcal{B}}^{u_n})$.
As observed in \cite{meh}, considering $V$ as a vector space over $F$, the map $q_{h,\mathcal{B}}$ is a system of quadratic forms over $F$.
If $h$ is even, we will consider $\mathcal{B}$ as a basis of $\symd_\lambda(D,\theta)$, instead of $\sym_\lambda(D,\theta)$.

\begin{prop}\label{st}
Let $(D,\theta)$ be a finite dimensional division algebra with involution of the first kind over a field $F$ and let $(V,h)$ be an even $\lambda$-hermitian space over $(D,\theta)$, where $\lambda=\pm1$.
Suppose that either $D\neq F$ or $\lambda\neq-1$.
If $h$ is regular then for every basis $\mathcal{B}$ of $\symd_\lambda(D,\theta)$ the system $q_{h,\mathcal{B}}$ is totally regular.
\end{prop}

\begin{proof}
In view of \cite[(2.6)]{knus}, the hypothesis $D\neq F$ or $\lambda\neq-1$ implies that $\symd_\lambda(D,\theta)\neq\{0\}$.
Hence, the result follows from \cite[(3.5)]{meh}.
\end{proof}

Let $(D,\theta)$ be a finite dimensional division algebra with involution of the first kind over a field $F$.
Let $K/F$ be a finite extension such that $D_K$ is a division ring.
Then for every $F$-basis $\mathcal{B}$ of $\symd_\lambda(D,\theta)$, the set $\{u\otimes1\mid u\in\mathcal{B}\}$ is a $K$-basis of
$\symd_\lambda((D,\theta)_K)=\symd_\lambda(D,\theta)\otimes K$.
We denote this basis by $\mathcal{B}_K$. 
\begin{thm}\label{main}
Let $(D,\theta)$ be a finite dimensional division algebra with involution of the first kind over a field $F$.
Let $K/F$ be a separable quadratic extension such that $D_K$ is a division ring.
Set $\lambda=\pm1$ and suppose that either $D\neq F$ or $\lambda\neq-1$.
For a  regular even $\lambda$-hermitian space $(V,h)$ over $(D,\theta)_K$ the following statements are equivalent.
\begin{itemize}
  \item [$(i)$] $h$ has a descent to $(D,\theta)$.
  \item [$(ii)$] There exists a basis $\mathcal{B}$ of $\symd_\lambda(D,\theta)$ for which $q_{h,\mathcal{B}_K}$ has a descent to $F$.
  \end{itemize}
\end{thm}

\begin{proof}
  If $h$ has a descent $h'$ to $(D,\theta)$, then for every basis $\mathcal{B}$ of $\symd_\lambda(D,\theta)$ we have $q_{h,\mathcal{B}_K}\simeq (q_{h',\mathcal{B}})_K$ by \cite[(3.4)]{meh}.
 This proves the implication $(i)\Rightarrow(ii)$.

Conversely, suppose that $(V,q_{h,\mathcal{B}_K})\simeq (V'_K,q'_K)$, where $\mathcal{B}=\{u_1,\cdots,u_m\}$ is a basis of $\symd_\lambda(D,\theta)$, $V'$ is a vector space over $F$ and $q'=(q'_1,\cdots,q'_m)$ is a quadratic map on $V'$.
We use induction on the right dimension $\dim_{D_K}V$ of $V$ over $D_K$. 
The assumption $D\neq F$ or $\lambda\neq-1$ implies that there exists $v\in V$ such that $h(v,v)\neq0$, i.e., $q_{h,\mathcal{B}_K}(v)\neq0$.
As $q_{h,\mathcal{B}_K}\simeq q'_K$ it follows that $q'(v_1)\neq0$ for some $v_1\in V'$.
Let $h_1:v_1D\times v_1D\rightarrow D$ be the $\lambda$-hermitian form satisfying
\[h_1(v_1,v_1)=\sum_{i=1}^m q_i'(v_1)u_i\in{\symd}_\lambda(D,\theta).\]
We may identify $(v_1\otimes1)D_K$ with a subspace of $V$, so that
\[h(v_1\otimes1,v_1\otimes1)=\sum_{i=1}^m q_{h,\mathcal{B}_K}^{u_i\otimes1}(v_1\otimes1)(u_i\otimes1)=\sum_{i=1}^m q_i'(v_1)(u_i\otimes1)=h_1(v_1,v_1)\otimes1.\]
Hence, $h|_{(v_1\otimes1)D_K\times (v_1\otimes1)D_K}\simeq(h_1)_K$.
If $\dim_{D_K}V=1$ then $h_1$ is a descent of $h$ and we are done. 
Otherwise, let
\[W=((v_1\otimes1)D_K)^{\perp_h}\subseteq V\quad {\rm and}\quad h'=h|_{W\times W}.\]
Since $h$ and $(h_1)_K$ are regular, we have $h\simeq (h_1)_K\perp h'$ by \cite[Ch. I, (3.6.2)]{knus1}, which implies that $q_{h,\mathcal{B}_K}\simeq q_{(h_1)_K,\mathcal{B}_K}\perp q_{h',\mathcal{B}_K}$.
Note that $q_{(h_1)_K,\mathcal{B}_K}\simeq (q_{h_1,\mathcal{B}})_K$, so $q_{(h_1)_K,\mathcal{B}_K}$ has a descent to $F$.
By Proposition \ref{st}, the systems $q_{h,\mathcal{B}_K}$, $q_{(h_1)_K,\mathcal{B}_K}$ and $q_{h',\mathcal{B}_K}$ are all totally regular.
It follows from Proposition \ref{desq} that $q_{h',\mathcal{B}_K}$ has a descent to $F$.
The pair $(W,h')$ is therefore an even $\lambda$-hermitian space over $(D,\theta)_K$ such that $q_{h',\mathcal{B}_K}$ has a descent to $F$.
By induction hypothesis, $h'$ has a descent to $(D,\theta)$.
Hence, $h\simeq (h_1)_K\perp h'$ has a descent to $(D,\theta)$, proving the result.
\end{proof}
\begin{rem}
The implication $(ii)\Rightarrow (i)$ in Theorem \ref{main} does not hold if $D=F$ and $\lambda=-1$.
Indeed, if $\car F\neq2$ then $\symd_{\lambda}(D,\theta)=\{0\}$ and the system $q_{h,\mathcal{B}_K}$ is trivial.
Hence, $q_{h,\mathcal{B}_K}$ does not give any information about $h$.

Suppose now that $\car F=2$.
Then $\sym_\lambda(D,\theta)=F$.
Write $K=F(\eta)$ for some $\eta\in K$ with $\eta^2+\eta=\delta\in F\setminus \wp(F)$, where $\wp(F)=\{x^2+x\mid x\in F\}$.
Suppose that the field $F$ satisfies $F\neq F^2$ and choose an element $a\in F\setminus F^2$.
Consider the diagonal bilinear form $h=\langle1,1+a\delta+a\eta\rangle$ over $K$ with the diagonal basis $\{v,w\}$, hence
\[h(v,v)=1,\quad h(w,w)=1+a\delta+a\eta \quad {\rm and}\quad h(v,w)=h(w,v)=0.\]
Then the determinant of $h$ is $(1+a+a^2\delta^2)F^{\times2}\in F^\times/F^{\times2}$.
Let $N_{K/F}:K\rightarrow F$ be the norm of $K$ over $F$.
Since
\[N_{K/F}(1+a\delta+a\eta)=1+a+a^2\delta^2\notin F^2,\]
$h$ has no descent to $F$.
However, taking $\mathcal{B}=\{1\}$, the system $q=q_{h,\mathcal{B}_K}$ is the quadratic form associated to $h$ (which is totally singular, i.e., its polar form is zero).
We have $q(v)=1$ and $q(\eta v+(1+\eta)w)=1+a\delta^2\in F$, implying that $q$ has a descent to $F$.
This proves the claim.
Note that this example also shows that the implication $(ii)\Rightarrow (i)$ in Theorem \ref{main} is not necessarily valid for non-even hermitian forms.
\end{rem}

\section{Hermitian forms over quaternion algebras}
Let $K/F$ be a quadratic separable extension with the nontrivial automorphism $\iota$.
For a central simple algebra $A$ over $K$, define the {\it conjugate algebra} $^\iota A=\{^\iota x\mid x\in A\}$ with operations
\[^\iota x+{^\iota y}= {^\iota}(x+y),\quad ^\iota x^\iota y={^\iota}(xy)\quad {\rm and}\quad ^\iota(\alpha x)=\iota(\alpha){^\iota x},\]
for $x,y\in A$ and $\alpha\in K$.
Let $s:{^\iota}A\otimes_KA\rightarrow{^\iota}A\otimes_KA$ be the $\iota$-semilinear map induced by $s({^\iota}x\otimes y)={^\iota} y\otimes x$ for $x,y\in A$.
The {\it corestriction} ${\co}_{K/F}(A)$ of $A$ is defined as
\[{\co}_{K/F}(A)=\{x\in{^\iota}A\otimes_KA\mid s(x)=x\}.\]
By \cite[(3.13 (4))]{knus}, $\co_{K/F}(A)$ is a central simple algebra over $F$.

\begin{lem}\label{quatdes}
Let $K/F$ be a separable quadratic extension and let $(Q,\sigma)$ be a quaternion algebra with involution of the first kind over $K$.
\begin{itemize}
  \item [$(i)$] If $\sigma$ is symplectic then $(Q,\sigma)$ has a descent to $F$ if and only if $\co_{K/F}(Q)$ splits.
  \item [$(ii)$] If $\sigma$ is orthogonal then $(Q,\sigma)$ has a descent to $F$ if and only if $\co_{K/F}(Q)$ splits and there exists $u\in\symd_{-1}(Q,\sigma)$ such that $u^2\in F$.
Furthermore, if these conditions are satisfied then there exists a descent $(Q',\sigma')$ of $(Q,\sigma)$ such that $u\in\symd_{-1}(Q',\sigma')$.
\end{itemize}
\end{lem}

\begin{proof}
$(i)$ According to \cite[(2.22)]{knus} and \cite[Ch. 8, (9.5)]{sch}, $Q$  has a descent to $F$ if and only if $\co_{K/F}(Q)$ splits.
By \cite[(2.21)]{knus}, the unique symplectic involution on $Q$ is the canonical involution.
Hence, if $\co_{K/F}(Q)$ splits then for every descent $Q'$ of $Q$, one has $(Q,\sigma)\simeq(Q',\gamma)_K$, where $\gamma$ is the canonical involution on $Q'$.
It follows that $(Q,\sigma)$ has a descent to $F$ if and only if $\co_{K/F}(Q)$ splits.

$(ii)$
If $\car F\neq2$, the claim follows from \cite[(2.4)]{dhedre} and its proof.
Otherwise, it can be found in \cite[(4.7)]{me}.
\end{proof}

Let $(Q,\sigma)$ be a quaternion division algebra with involution of type $\varepsilon$ over a field $F$ and set $\lambda=-\varepsilon$.
Then $\dim_F\symd_\lambda(Q,\sigma)=1$ by \cite[(2.6)]{knus}.
Choose a basis $\mathcal{B}=\{u\}$ of $\symd_\lambda(Q,\sigma)$.
Let $(V,h)$ be an even $\lambda$-hermitian space on $(Q,\sigma)$.
Then the quadratic map $q_{h,\mathcal{B}}$ reduces to a quadratic form, which we denote it by
$q_{h,u}$.
Note that the form $q_{h,u}$ is uniquely determined, up to a scalar factor.
Also, if $\sigma$ is symplectic then $\symd_1(Q,\sigma)=F$.
Taking $u=1$, the form $q_{h,u}$ coincides with the Jacobson's trace form of $h$, introduced in \cite{jacob}.
We will simply denote this from by $q_h$.

Let $K/F$ be a quadratic separable extension and let $(D,\theta)$ be a division algebra with involution of the first kind over $K$.
If $(D,\theta)$ has a descent $(D',\theta')$ to $F$ then Theorem \ref{main} gives a criterion for hermitian spaces over $(D,\theta)$ to have a descent to $(D',\theta')$.
In the case where $(D,\theta)$ is a quaternion division algebra with involution of type $\varepsilon$ over $K$ and $\lambda=-\varepsilon$,
using Lemma \ref{quatdes}, we state the following criterion for $\lambda$-hermitian spaces over $(D,\theta)$ to have a descent to $F$.

\begin{thm}\label{main2}
Let $K/F$ be a separable quadratic extension and let $(Q,\sigma)$ be a quaternion division algebra with involution of type $\varepsilon$ over $K$.
Set $\lambda=-\varepsilon$.
Let $(V,h)$ be a regular $\lambda$-hermitian space over $(Q,\sigma)$.
\begin{itemize}
  \item [$(i)$] If $\sigma$ is symplectic then $h$ has a descent to $F$ if and only if $\co_{K/F}(Q)$ splits and $q_h$ has a descent to $F$.
  \item [$(ii)$] If $\sigma$ is orthogonal then $h$ has a descent to $F$ if and only if $\co_{K/F}(Q)$ splits and $q_{h,u}$ has a descent to $F$ for some $u\in\symd_{-1}(Q,\sigma)$ with $u^2\in F$.
\end{itemize}
  \end{thm}

\begin{proof}
$(i)$ If $h$ has a descent to $F$ then there exists a quaternion $F$-algebra with involution $(Q',\sigma')$ and a hermitian form $h'$ on $(Q',\sigma')$ such that $(Q,\sigma)\simeq(Q',\sigma')_K$ and $h\simeq h'_K$.
By Lemma \ref{quatdes} $(i)$, $\co_{K/F}(Q)$ splits.
We also have $q_h\simeq (q_{h'})_K$, proving the `only if' implication.

Conversely, if $\co_{K/F}(Q)$ splits then $(Q,\sigma)$ has a descent $(Q',\sigma')$ to $F$, thanks to Lemma \ref{quatdes} $(i)$.
Since $q_h$ has a descent to $F$, the result follows from Theorem \ref{main} by taking $\mathcal{B}=\{1\}$.

$(ii)$ Suppose that $h$ has a descent to $F$.
Then there exists a quaternion $F$-algebra with involution $(Q',\sigma')$ and a skew hermitian form $h'$ on $(Q',\sigma')$ such that $(Q,\sigma)\simeq(Q',\sigma')_K$ and $h\simeq h'_K$.
By \cite[(3.13 (5))]{knus}, $\co_{K/F}(Q)$ splits.
Let $u'\in\symd_{-1}(Q',\sigma')$ be a nonzero element.
Considering the isomorphism $(Q,\sigma)\simeq(Q',\sigma')_K$ as an identification, we have $u:=u'\otimes1\in\symd_{-1}(Q,\sigma)$, $u^2\in F$ and $q_{h,u}\simeq (q_{h',u'})_K$.

Conversely, suppose that $\co_{K/F}(Q)$ splits and $(V,q_{h,u})\simeq (V',q')_K$, where $u\in\symd_{-1}(Q,\sigma)$ is an element with $u^2\in F$ and $(V',q')$ is a quadratic space over $F$.
By Lemma \ref{quatdes} $(ii)$, $(Q,\sigma)$ has a descent $(Q',\sigma')$ to $F$ with $u\in\symd_{-1}(Q',\sigma')$.
As $\symd_{-1}(Q,\sigma)$ is one-dimensional, the result follows from Theorem \ref{main} by taking $\mathcal{B}=\{u\}$.
\end{proof}

\noindent{\sc A.-H. Nokhodkar, {\tt
    a.nokhodkar@kashanu.ac.ir},\\
Department of Pure Mathematics, Faculty of Science, University of Kashan, P.~O. Box 87317-51167, Kashan, Iran.}

\end{document}


\maketitle

\begin{abstract}
\noindent
There exists a mistake in the proof of Theorem 4.2.
We present a new proof of this theorem, which shows that the main results of the paper are still true.\\
\\
\end{abstract}

\section{Introduction}
In \cite[\S3]{mem}, it was claimed that for a finite field extension $K/F$, a system of quadratic forms $q=(q_1,\cdots, q_m)$ over $K$ has a descent to $F$ if and only if every $q_i$ has a descent to $F$.
This is incorrect in general.
For example, let $K/F$ be a separable quadratic extension and let $\alpha\in K^2\setminus F$.
The $2$-fold system of quadratic forms $(\langle1\rangle,\langle\alpha\rangle)$ over $K$ has no descent to $F$, while both forms $\langle1\rangle$ and $\langle\alpha\rangle$ can be descended to $F$.
This claim is used in the proof of Proposition 3.2, which in turn plays a critical role in the proof of the implication $(ii)\Rightarrow(i)$ of Theorem 4.2.
Hence, the given proof of Theorem 4.2 is not true.
We present a proof of the implication $(ii)\Rightarrow(i)$ of Theorem 4.2 with a different approach, which does not depend on Proposition 3.2.
Since Proposition 3.2 was only used in the proof of Theorem 4.2, one can remove this proposition, with no effect on the other results of the paper.

\section{The correction}
Throughout this note, $(D,\theta)$ is a finite dimensional division algebra with involution of the first kind over a field $F$ and $\lambda=\pm1$.
We also fix $K/F$ as a separable quadratic extension such that $D_K$ is a division ring.

Let $q=(q_1,\cdots,q_n)$ be a system of quadratic forms over $K$ and let $s:K\rightarrow F$ be a nonzero $F$-linear map.
We call the system of quadratic forms $s_*(q):=(s_*(q_1),\cdots,s_*(q_n))$ over $F$ the {\it transfer} of $q$.
Similarly, one can define the {\it  transfer} of a hermitian form as follows: let $(V,h)$ be a $\lambda$-hermitian space over $(D,\theta)_K$.
Extend $s$ to a $D$-linear map $s_D:D_K\rightarrow D$ and define $s_*(h):V\times V\rightarrow D$ via
\[s_*(h)(u,v)=s_D(h(u,v))\quad {\rm for \ all \ } u,v\in V.\]
Clearly, $s_*(h)$ is a $\lambda$-hermitian form over $(D,\theta)$.
It is easy to see that if $h$ is regular then so is $s_*(h)$.
It is also easily verified that if $h$ is metabolic then so is $s_*(h)$ (see \cite[p. 362]{bayer}).

\begin{lem}\label{lq}
Let $(V,h)$ be a $\lambda$-hermitian form over $(D,\theta)_K$ and let $s:K\rightarrow F$ be a nonzero $F$-linear map.
Then $s_*(q_{h,\mathcal{B}_K})=q_{s_*(h),\mathcal{B}}$ for every basis $\mathcal{B}$ of $\sym_\lambda(D,\theta)$.
\end{lem}

\begin{proof}
Write $\mathcal{B}=\{u_1,\cdots,u_m\}$ for some $u_1,\cdots,u_m\in \sym_\lambda(D,\theta)$.
It suffices to show that $s_*(q_{h,\mathcal{B}_K}^{u_i\otimes1})=q_{s_*(h),\mathcal{B}}^{u_i}$ for all $i=1,\cdots,m$.
Fix an index $1\leqslant i\leqslant m$ and a vector $v\in V$.
Write $K=F(\eta)$ for some $\eta\in K$.
Since $\sym_\lambda((D,\theta)_K)=\sym_\lambda(D,\theta)\otimes K$, identifying $D$ with a subalgebra of $D_K$, one can write $h(v,v)=\alpha+\beta\eta$, where $\alpha,\beta\in \sym_\lambda(D,\theta)$.
Write
\[\textstyle\alpha=\sum_{j=1}^ma_ju_j \quad {\rm and} \quad \beta=\sum_{j=1}^mb_ju_j,\]
for some $a_j,b_j\in F$, $j=1,\cdots,m$.
Extend $\pi_i:\sym_\lambda(D,\theta)\rightarrow F$ to a $K$-linear map $(\pi_i)_K:\sym_\lambda((D,\theta)_K)\rightarrow K$.
Then
\begin{align*}
\textstyle
s_*(q_{h,\mathcal{B}_K}^{u_i\otimes1})(v)&=s((\pi_i)_K(h(v,v)))\\
&=\textstyle s((\pi_i)_K(\sum_{j=1}^m(a_ju_j+b_ju_j\eta)))\\
&=s(a_i+b_i\eta)=a_is(1)+b_is(\eta)\\
&\textstyle=\pi_i(\sum_{j=1}^ma_ju_js(1)+\sum_{j=1}^mb_ju_js(\eta))\\
&=\pi_i(\alpha s(1)+\beta s(\eta))=\pi_i(s_D(\alpha+\beta\eta))\\
&=\pi_i(s_D(h(v,v)))=\pi_i(s_*(h)(v,v))=q_{s_*(h),\mathcal{B}}^{u_i}(v).\qedhere
\end{align*}
\end{proof}
We recall that a regular even $\lambda$-hermitian form is metabolic if and only if it is hyperbolic (see \cite[Ch. I, (3.7.3)]{knus1}).
\begin{lem}\label{lhyp}
Let $(V,h)$ be a regular even $\lambda$-hermitian form over $(D,\theta)_K$ and let $s:K\rightarrow F$ be a nonzero $F$-linear map with $s(1) = 0$.
Then $h$ has a descent to $(D,\theta)$ if and only if $s_*(h)$ is hyperbolic.
\end{lem}
\begin{proof}
Observe first that $s_*(h)$ is even and regular.
Suppose that $(V,h)\simeq(V',h')_K$ for some $\lambda$-hermitian space $(V',h')$ over $(D,\theta)$.
We may identify $V'\subseteq V$ and $D\subseteq D_K$.
The assumption $s(1)=0$ then implies  that $s_D|_D=0$, hence, $s_*(h)|_{V'\times V'}=0$.
Since $\dim_DV'=\frac{1}{2}\dim_DV$, one concludes that  $s_*(h)$ is metabolic and therefore hyperbolic.

To prove the converse, write $h\simeq h_{1}\perp h_2$, where $h_{1}$ is anisotropic and $h_2$ is hyperbolic.
Clearly, $h_2$ has a descent to $(D,\theta)$.
Hence, it suffices to show that  $h_1$ has a descent to $(D,\theta)$.
Since $h_2$ is metabolic, $s_*(h_2)$ is metabolic and therefore hyperbolic.
It follows that $s_*(h_{1})$ is also hyperbolic, because $s_*(h)\simeq s_*(h_{1})\perp s_*(h_2)$.
Since $h_{1}$ is anisotropic and the map $s_D$ is nonzero with $s_D|_D= 0$, one concludes that $h_1$ represents a nonzero element $\alpha\in D$.
By \cite[Ch. I, (3.6.2)]{knus1}, there exists a regular (even) $\lambda$-hermitian form $h'$ over $(D,\theta)_K$ such that  $h_{1}\simeq\langle \alpha\rangle_K \perp h'$.
Since $\langle\alpha\rangle_K$ has a descent to $(D,\theta)$, the above argument shows that $s_*(\langle \alpha\rangle_K)$ is hyperbolic, hence $s_*(h')$ is also hyperbolic.
The result now follows by induction on $\dim_Dh_1$.
\end{proof}

As in \cite[Ch. 9]{pfister}, we say that a system $(V,q)$ of quadratic forms over $F$ is {\it metabolic}
if there exists a subspace $L$ of $V$ with $\dim_F L \geqslant\frac{1}{2}\dim_F V$ such that $q|_L = 0$.

\begin{lem}\label{31}
Let $(V,q)$ be a system of quadratic forms over $K$ and let $s:K\rightarrow F$ be a nonzero $F$-linear map with $s(1) = 0$.
If $q$ has a descent to $F$ then $s_*(q)$ is metabolic.
\end{lem}
\begin{proof}
Write $(V,q)\simeq(V',q')_K$, where $(V',q')$ is  a system of quadratic forms over $F$.
Considering $V'$ as an $F$-subspace of $V$, one has $s_*(q)|_{V'}=0$.
Since $\dim_FV'=\frac{1}{2}\dim_FV$, the system $s_*(q)$ is metabolic.
\end{proof}

We are now able to correct the proof of Theorem 4.2 of \cite{mem} as follows:

\noindent
{\it Proof of implication $(ii)\Rightarrow (i)$ of \cite[(4.2)]{mem}.}
Let $s:K\rightarrow F$ be a nonzero $F$-linear map with $s(1) = 0$.
If $q_{h,\mathcal{B}_K}$ has a descent to $F$ then
$s_*(q_{h,\mathcal{B}_K})$ is metabolic by \cref{31}.
Hence, $q_{s_*(h),\mathcal{B}}$ is also metabolic by \cref{lq}.
The assumption  $D \neq F$ or $\lambda \neq -1$ implies that either $\car F\neq2$ or $D\neq F$.
It follows from \cite[(4.2)]{meh} that $s_*(h)$ is metabolic and therefore hyperbolic.
The result now follows from \cref{lhyp}.
\\
\\
{\bf Acknowledgements.} The author is grateful to Uriya First, who pointed out the mistake in the proof of Proposition 3.2.
This research is partially supported by the University of Kashan under the grant number 890193/2.

\noindent{\sc A.-H. Nokhodkar, {\tt
    a.nokhodkar@kashanu.ac.ir},\\
Department of Pure Mathematics, Faculty of Science, University of Kashan, P.~O. Box 87317-51167, Kashan, Iran.}